\documentclass[12pt]{article}
\usepackage{latexsym}
\usepackage{amsmath, amsthm}
\usepackage[all,dvips,ps]{xy}
\usepackage{amsfonts}
\usepackage{mathrsfs}
\usepackage{tikz}
\newtheorem{prop}{Proposition}[section]
\newtheorem{lem}{Lemma}[section]
\newtheorem{thm}{Theorem}[section]
\newtheorem{cor}{Corollary}[section]

\newtheorem{exa}{Example}[section]

\newtheorem{defn}{Definition}[section]

\newcommand{\tra}{\mathrm{trace}\,}
\newcommand{\diag}{\mathrm{diag}\,}
\newcommand{\rank}{\mathrm{rank}\,}

\newcommand{\Rs}{\mathbb{R}}

\newcommand{\dE}{\mathrm{dim}_E(G)}
\newcommand{\dS}{\mathrm{dim}_S(G)}
\newcommand{\dJ}{\mathrm{dim}_J(G)}
\newcommand{\mmin}{\mu_{\mathrm{min}} }
\newcommand{\mmax}{\mu_{\mathrm{max}} }
\newcommand{\mminb}{\bar{\mu}_{\mathrm{min}} }
\newcommand{\mmaxb}{\bar{\mu}_{\mathrm{max}} }

\newcommand{\bz}{{\bf 0} }

\newcommand{\bpr}{{\bf Proof.} \hspace{1 em}}
\newcommand{\epr}{ \\ \hspace*{4.5in} $\Box$ }
\newcommand{\beq}{ \begin{equation} }
\newcommand{\eeq}{ \end{equation} }
\newcommand{\bt}{ \begin{tabular} }
\newcommand{\et}{ \end{tabular} }

\begin{document}

\bibliographystyle{plain}
\title{On Representations of Graphs as Two-Distance Sets}
\vspace{0.3in}
        \author{ A. Y. Alfakih
  \thanks{E-mail: alfakih@uwindsor.ca}
  \\
          Department of Mathematics and Statistics \\
          University of Windsor \\
          Windsor, Ontario N9B 3P4 \\
          Canada
}

\date{ \today}
\maketitle

\noindent {\bf AMS classification: 05C50, 05C62, 51K05, 15B48.}

\noindent {\bf Keywords:}Two-distance sets, Euclidean, spherical and
$J$-spherical representations, Euclidean distance matrices (EDMs), spherical EDMs.
\vspace{0.1in}

\begin{abstract}
Let $\alpha \neq \beta$ be two positive scalars.
A Euclidean representation of a simple graph $G$ in $\Rs^r$ is a mapping of the nodes of $G$ into points in $\Rs^r$
such that the squared Euclidean distance between any two points  is $\alpha$ if the corresponding nodes are adjacent
and $\beta$ otherwise. A Euclidean representation is spherical if the points lie on an $(r-1)$-sphere,
and is $J$-spherical if this sphere has radius 1 and $\alpha=2 < \beta$.
Let $\dE$, $\dS$ and $\dJ$ denote, respectively, the smallest dimension $r$ for which
$G$ admits a  Euclidean, spherical and $J$-spherical representation.

In this paper, we extend and simplify the results of Roy \cite{roy10} and Nozaki and Shinohara \cite{ns12} by
deriving exact simple formulas for $\dE$ and $\dS$ in terms of the eigenvalues of $V^TAV$, where $A$ is the adjacency matrix
of $G$ and $V$ is the matrix whose columns form an orthonormal basis for the orthogonal complement of the vector of all 1's.
We also extend and simplify the results of Musin \cite{mus18} by deriving explicit formulas for determining the
$J$-spherical representation of $G$ and for determining $\dJ$ in terms of the largest eigenvalue of $\bar{A}$,
the adjacency matrix of the complement graph $\bar{G}$.
As a by-product, we obtain several other related results and in particular we answer a question raised by Musin in \cite{mus18}.
\end{abstract}

\section{Introduction}

Let $G$ be a simple graph, i.e., no loops and no multiple edges, on $n$ nodes. A {\em Euclidean representation} of $G$ in
 $\Rs^r $, the $r$-dimensional Euclidean space,
is an $n$-point configuration $p^1, \ldots, p^n$ in $\Rs^r$ such that: for all $i,j=1,\ldots, n$, we have
\[
|| p^i - p^j ||^2 = \left\{ \begin{array}{ll} \alpha & \mbox{ if } \{i,j\} \in E(G), \\
                                               \beta  & \mbox{ if } \{i,j\} \not \in E(G), \end{array} \right.
\]
for two distinct positive scalars $\alpha$ and $\beta$, where $||x||^2 = x^T x $ and $E(G)$ is the set of edges of $G$.
In other words, the points $p^1, \ldots, p^n$ form a two-distance set.
The {\em Euclidean representation number}  \cite{roy10} of $G$, denoted by $\dE$,
is the smallest $r$ for which $G$ admits a Euclidean representation in $\Rs^r$.
A Euclidean representation of $G$ in $\Rs^r$  is said to be a {\em spherical representation} of $G$ in $\Rs^r$ if the points
$p^1, \ldots, p^n$ lie on an $(r-1)$-sphere in $\Rs^r$. Moreover,
the {\em spherical representation number}  of $G$, denoted by $\dS$,
is the smallest $r$ for which $G$ admits a spherical representation in $\Rs^r$.
In the special case of a spherical representation of $G$, where the sphere has unit radius and $\alpha = 2 < \beta$,
the spherical representation  is said to be a {\em $J$-spherical representation} \cite{mus18}.
In the same manner, the {\em $J$-spherical representation number}  of $G$, denoted by $\dJ$,
is the smallest $r$ for which $G$ admits a $J$-spherical representation in $\Rs^r$.
Evidently
\[
\dE \leq \dS \leq \dJ.
\]

Einhorn and Schoenberg \cite{es66a} gave exact formulas for $\dE$ in terms of the multiplicities of certain roots
of the discriminating polynomial defined in (\ref{dispoly}).
They \cite{es66b} also determined all two-distance sets in dimensions two and three.
A full classification of all maximal two-distance sets in dimension $r$ for all $r \leq 7$ is given in \cite{lis97}.
Recently, there has been a renewed interest in the problems of
determining $\dE$, $\dS$ and $\dJ$  \cite{roy10,ns12,mus18}.
Roy \cite{roy10} derived bounds on $\dE$ using the multiplicities of the smallest and the second  smallest distinct eigenvalues
of $A$, the adjacency matrix of $G$. He also gave exact formulas for $\dE$ using the main angles of the graph.
Nozaki and Shinohara \cite{ns12} considered the problem of determining $\dS$ and,
using Roy's results, they obtained necessary and sufficient conditions for a Euclidean representation of $G$ to be spherical.
Musin \cite{mus18} considered the problem of determining $\dJ$ and proved that
any graph which is neither complete nor null admits a unique, up to an isometry, $J$-spherical representation.
He also obtained exact formulas for $\dS$ and $\dJ$ in terms of the multiplicities of the roots of
a polynomial defined by the Cayley-Menger determinant.
Finally, we should point out that a classification of all two-distance sets in dimension four is given in \cite{szo18}.

In this paper, we extend and simplify the results of Roy \cite{roy10} and Nozaki and Shinohara \cite{ns12} by
deriving exact simple formulas for $\dE$ and $\dS$ in terms of the multiplicities of the
smallest and the largest eigenvalues of the $(n-1) \times (n-1)$ matrix
$V^TAV$, where $A$ is the adjacency matrix
of $G$ and $V$, defined in (\ref{defV}), is the matrix whose columns form an orthonormal basis for the orthogonal
complement of the vector of all 1's. This is made possible by using projected Gram matrices for representing
$n$-point configurations.  As a by-product, we obtain a characterization of $(0-1)$ Euclidean distance matrices (Theorem \ref{thmcmpar}).

We also extend and simplify the results of Musin \cite{mus18} by deriving explicit formulas for determining the
$J$-spherical representation of $G$ and for determining $\dJ$ in terms of the largest eigenvalue of $\bar{A}$,
the adjacency matrix of the complement graph $\bar{G}$, and
its multiplicity. This is made possible by the extensive use of the theory of Euclidean distance matrices.
We also answer a question raised by Musin in \cite{mus18}.

The remainder of this paper is organized as follows. Section 2 presents the background we need from
Euclidean matrices (EDMs), spherical EDMs, projected Gram matrices
 and Gale transform. In Section 3, we present some of the spectral properties of the
matrix $V^TAV$ since this matrix plays a key role in determining $\dE$ and $\dS$. Sections 4, 5 and 6  discuss,
respectively, Euclidean, spherical and $J$-spherical representation of graph $G$.

\subsection{Notation}
We collect here the notation used throughout  the paper.
$e_n$ and $E_n$ denote, respectively, the vector of all 1's  in $\Rs^n$ and the $n \times n$ matrix of all 1's.
The subscript is omitted if the dimension is clear from the context.
The identity matrix of order $n$ is denoted by $I_n$.
The zero matrix or zero vector of appropriate dimension is denoted  by $\bz$.
For a matrix $A$, $\diag(A)$ denotes the vector consisting of the diagonal entries of $A$.
$m(\lambda)$ denotes the multiplicity of eigenvalue $\lambda$.

$K_n$ denotes the complete graph on $n$ nodes.
The adjacency matrix of a graph $G$ is denoted by $A$, and the adjacency matrix of the complement graph $\bar{G}$
is denoted by $\bar{A}$. $\mmin$ and $\mmax$ denote, respectively, the minimum and the maximum eigenvalues of $V^TAV$.
Likewise, $\mminb$ and $\mmaxb$ denote, respectively, the minimum and the maximum eigenvalues of $V^T\bar{A}V$.
Finally, PSD and PD stand for positive semidefinite and positive definite.

\section{Preliminaries}

The theory of Euclidean distance matrices (EDMs) provides a natural and powerful tool for determining $\dE$, $\dS$ and $\dJ$.
In this section, we present the necessary background concerning  EDMs, spherical EDMs, projected Gram matrices
and Gale matrices. For a comprehensive treatment of these topics and  EDMs in general, see the monograph \cite{alf18}.

\subsection{EDMs}

An $n \times n$ matrix $D=(d_{ij})$ is said to be an EDM if there exist points $p^1,\ldots,p^n$ in some Euclidean
space such that
\[
d_{ij}= || p^i - p^j ||^2 \mbox{ for all } i,j=1,\ldots, n,
\]
$p^1,\ldots,p^n$ are called the {\em generating points} of $D$ and the dimension of their affine span is
called the {\em embedding dimension} of $D$. Let $D$ be an EDM of embedding dimension $r$. We always assume
throughout this paper that the generating points of $D$ are in $\Rs^r$. Hence, the $n \times r$ matrix
\[
P= \left[ \begin{array}{c} (p^1)^T \\ \vdots \\ (p^n)^T \end{array} \right]
\]
has full column rank. $P$ is called a {\em configuration matrix} of $D$.

Let $e$ be the vector of all 1's in $\Rs^n$ and let $E=ee^T$.
The following theorem is a well-known characterization of EDMs \cite{sch35,yh38, gow85,cri88}.
\begin{thm}  \label{thmbasic}
Let $D$ be an $n \times n$ real symmetric matrix whose diagonal entries are all 0's
and let $s \in \Rs^n$ such that $e^Ts= 1$.
Then $D$ is an EDM
if and only if $D$ is negative semidefinite on $e^{\perp}$, the orthogonal complement of $e$ in $\Rs^n$;
i.e., iff
\beq \label{defs}
B= -\frac{1}{2} (I - e s^T) D (I - s e^T )
\eeq
is positive semidefinite (PSD),
in which case, the embedding dimension of $D$ is given by $\rank(B)$.
\end{thm}

Note that $B$, which can be factorized as $B=PP^T$, is the Gram matrix of the generating points of $D$,
or the Gram matrix of $D$ for short. Moreover, $Bs=\bz$ and hence $P^T s=\bz$.
It is well known \cite{gow85} that if $D$ is a nonzero EDM, then $e$ lies in the column space of $D$. Hence, there exists $w$ such that
\beq  \label{defw}
D w = e.
\eeq
Two choices of vector $s$ in (\ref{defs}) are of particular interest to us. First, $s=e/n$. This choice fixes the origin at the centroid of
the generating points of $D$ and thus the corresponding Gram matrix satisfies $Be=\bz$. Second, $s=2 w$, where $w$ is as defined in (\ref{defw}).
This choice, as we will see in Section 6, is particularly useful when the radius of a spherical EDM is known.

Assume that $B$, the Gram matrix of $D$, satisfies $Be= \bz$.
Let $V$ be the $n \times (n-1)$ matrix whose columns form an orthonormal basis of $e^{\perp}$; i.e.,
$V$ satisfies
\beq \label{defV}
V^T e = \bz \mbox{ and } V^T V = I_{n-1}.
\eeq
Hence, $VV^T = I_n - E/n$ is the orthogonal projection on  $e^{\perp}$. Thus,
$-2B = V V^T D V V^T$. Let
\beq \label{defX}
X = V^T B V = - \frac{1}{2} V^T D V,
\eeq
and thus, $B = V X V^T$. Moreover,  it readily follows that $B$ is PSD of rank $r$ iff $X$ is PSD of rank $r$.
As a result, $X$ is called the {\em projected Gram matrix} of $D$. Consequently, a real symmetric matrix with $\diag(D)=\bz$
is an EDM of embedding dimension $r$ if and only if its projected Gram matrix $X$ is PSD of rank $r$.

It should be pointed out that $V$ as defined in (\ref{defV}) is not unique.
One such choice of $V$ is
\[
 V = \left[ \begin{array}{c} y e^T_{n-1}   \\  I_{n-1}+ x E_{n-1} \end{array} \right], \mbox{ where } y=\frac{-1}{\sqrt{n}} \mbox{ and } x=\frac{-1}{n+\sqrt{n}}.
\]
Another choice of $V$, which we use in the sequel and is particularly convenient
when dealing with block matrices, is
 \beq  \label{defVs}
 V= \left[ \begin{array}{ccc} V'_3 & \bz  & a e_3 \\ \bz  & V'_{n-3} & b e_{n-3} \end{array} \right].
\eeq
Here, $V'_3$ and $V'_{n-3}$ are, respectively, $3 \times 2$ and $(n-3) \times (n-4)$ matrices satisfying (\ref{defV}),
\[
a= \left( \frac{n-3}{3n}\right)^{1/2} \;\;\; , \;\;\;  b= -  \left( \frac{3}{n(n-3)}\right)^{1/2}
\]
and
\[
V'_3 = \left[ \begin{array}{cc} y & y  \\ 1+x & x \\ x & 1+x \end{array} \right], \;\; y=\frac{-1}{\sqrt{3}} \mbox{ and } x=\frac{-1}{3+\sqrt{3}}.
\]

Gale transform \cite{gal56, gru67}, or Gale matrix, plays an important role in theory of EDMs.  Let $Z$ be the $n \times (n-r-1)$ matrix
whose columns form a basis of the null space of
\[
\left[ \begin{array}{c} P^T \\ e^T \end{array} \right],
\]
where $P$ is a configuration matrix of $D$.
Then $Z$ is called a {\em Gale matrix} of $D$. The following lemma establishes the relationship between Gale matrix $Z$
and the null space of the projected Gram matrix $X$.

\begin{lem}[\cite{alf01}] \label{lemZVU}
Let $D$ be an $n \times n$ EDM of embedding dimension $r \leq n-2$ and let $X$ be the projected Gram matrix of $D$.
Further, let $U$ be the matrix whose columns form an orthonormal basis of the null space of $X$.  Then $VU$ is a Gale matrix of $D$,
where $V$ is as defined in (\ref{defV}).
\end{lem}

\subsection{Spherical EDMs}

An EDM $D$ is said to be {\em spherical} if its generating points lie on a sphere. We denote the radius of the generating points of
a spherical EDM  $D$ by $\rho$ and we will refer to it as the radius of $D$. Among the many different characterizations of spherical EDMs,
the ones that are relevant to this paper are given in the following theorem.

\begin{thm} \label{thmspher}
Let $D$ be an $n \times n$ EDM of embedding dimension $r$ and let $Dw =e$. Let $P$ and $Z$ be, respectively, a configuration matrix
and a Gale matrix of $D$ and assume that $P^Te=\bz$.
If $r = n-1$, then $D$ is  spherical. Otherwise, if $r \leq n-2$, then
the following statements are equivalent:
\begin{enumerate}
\item $D$ is spherical,
\item $DZ = \bz$.
\item $\rank(D)=r+1$.
\item there exists $a \in \Rs^r$ such that
\[
P a = \frac{1}{2} (I - \frac{E}{n}) \diag(PP^T)
\]
in which case, the generating points of $D$ lie on a sphere centered at $a$ and with radius
\[
 \rho = \left(a^Ta + \frac{e^T D e}{2 n^2} \right)^{1/2}.
\]
\item $e^Tw >0$, in which case, the radius of $D$ is given by
\beq
\rho = \left(\frac{1}{2 e^T w } \right)^{1/2}.
\eeq
\end{enumerate}
\end{thm}
The equivalence between Statement 1 and 2 was proven by Alfakih and Wolkowicz \cite{aw02}.
The equivalence between Statement 1 and 3 was proven by Gower \cite{gow85}.
The equivalence between Statement 1 and 4 was proven by Tarazaga et al \cite{thw96}.
Finally, the equivalence between Statement 1 and 5 was proven by Gower \cite{gow82,gow85}.

An interesting subclass of spherical EDMs is that of regular EDMs. A spherical EDM $D$ is said to
be {\em regular} if the center of the sphere containing the generating points of $D$ coincides with
the centroid of these points. Regular EDMs are characterized \cite{ht93} as those EDMs which have $e$ as an eigenvector.
It is easy to see that an $n \times n$ regular EDM $D$ has radius
\[
\rho = \left(\frac{e^T D e}{2n^2} \right)^{1/2}.
\]

\section{Spectral Properties of $V^TAV$}

Cluster graphs and complete multipartite graphs play a special role in this paper.
Graph $G$ is said to be a {\em cluster graph} if it is the disjoint union of complete graphs. Note that
$K_1$, the graph consisting of a single isolated node, is considered complete.
The complement of a cluster graph $G$  is called a {\em complete multipartite graph}.
Thus, the vertices of a complete multipartite graph can be partitioned into independent sets.
We often denote a complete multipartite graph by $K_{n_1,\ldots,n_s}$, where
$n_1, \ldots, n_s$ are the sizes of its independent sets.
Let $P_3$, or $K_{1,2}$, denote the graph consisting of a path on 3 nodes. Then
it is well known that $G$ is a cluster graph if and only if it is $P_3$-free, i.e., it
has no $P_3$ as an induced subgraph. As a result,
$G$ is a complete multipartite graph if and only if it is  $\overline{P_3}$-free.
It should be pointed out that $K_n$ is both a cluster graph and a complete multipartite graph, and thus
the null graph, $\overline{K_n}$, is also both a cluster graph and a complete multipartite graph.

Let $A$ denote the adjacency matrix of $G$ and let
$\mmin$ and $\mmax$ denote, respectively, the minimum and the maximum eigenvalues of $V^TAV$.
Note that if $A \neq \bz$, then $\mmin < 0$ since $\tra(V^TAV) = -e^TAe/n$.

\begin{prop} \label{propinterl}
Let $\lambda_n(A)$ and $\lambda_{n-1}(A)$ denote the smallest and the second smallest eigenvalues of $A$.
Also, let $\lambda_1(A)$ and $\lambda_{2}(A)$ denote the largest and the second largest eigenvalues of $A$.
 Then
\[
  \lambda_{1} \geq \mmax \geq \lambda_2(A) \; \mbox{ and } \; \lambda_{n-1} \geq \mmin \geq \lambda_n(A).
\]
\end{prop}

\bpr Let $Q= [ e/\sqrt{n} \;\; V]$. Then
\[
  Q^T A Q = \left[ \begin{array}{cc} e^T A e/n & e^T A V/ \sqrt{n} \\ V^TA e /\sqrt{n} & V^TAV \end{array} \right].
\]
The result follows from the interlacing theorem since $Q$ is orthogonal.
\epr

The following proposition is an immediate consequence of the proof of Proposition \ref{propinterl}.

\begin{prop}\label{propreg}
Let $A$ denote the adjacency matrix of a $k$-regular graph. Then the eigenvalues of $A$ are exactly those
of $V^TAV$ in addition to the eigenvalue $k$.
\end{prop}

\begin{lem} \label{lemmaxE-I}
Let $A$ denote the adjacency matrix of graph $G$.
Then
\begin{enumerate}
\item $\mmax = -1$ if and only if $A = E-I$.
\item $\mmax = 0$ if and only if $A$ is an EDM and $A \neq E-I$.
\end{enumerate}
\end{lem}

\bpr
Clearly, $\mmax \leq 0 $ iff ($-V^T A V $) is PSD iff $A$ is an EDM.
Now if $A=E-I$, then obviously $\mmax = -1$. On the other hand, if $\mmax = -1$,
then ($-V^TAV$) is  positive definite (PD) and thus
$A$ is an EDM of embedding dimension $n-1$.
Assume, by way of contradiction, that $A \neq E-I$. Then at least one off-diagonal entry of $A$ is zero,
 and thus at least two of
the generating points of $A$ coincide. Accordingly, the embedding dimension of $A$ is $\leq n-2$,
a contradiction. Therefore, $A=E-I$.
Also, we conclude that if $A$ is an EDM and $A \neq E-I$, then
the embedding dimension of $A$ is $\leq n-2$,
i.e., rank($V^T A V) \leq n-2$ and hence $\mmax = 0$ since ($-V^TAV$) is PSD.
\epr

The following theorem is a characterization of ($0-1$) EDMs.

\begin{thm} \label{thmcmpar}
Let $A$ denote the adjacency matrix of graph $G$. Then $A$ is an EDM
if and only if $G$ is a complete multipartite graph.
\end{thm}

\bpr
Assume that $G$ is a complete multipartite graph and
assume that the nodes of $G$ are partitioned into $s$ independent sets. Then obviously $A$
is an EDM whose generating points have the property that $p^i=p^j$ if and only if
nodes $i$ and $j$ belong to the same independent set.

To prove the other direction, assume that $G$ is not a complete multipartite graph. Then $G$ has
$\overline{P_3}$ as an induced subgraph. Wlog assume that the nodes of $P_3$ are 1, 2, and 3. Therefore,
 the third leading principal submatrix of $A$ is
\[
 \left[ \begin{array}{ccc} 0 & 0 & 1 \\ 0 & 0 & 0 \\ 1 & 0 & 0 \end{array} \right].
\]
Let $V$ be as defined in (\ref{defVs}). Then
the second leading  principal submatrix of $V^TAV$ is
\[
{V'}_3^T  \left[ \begin{array}{ccc} 0 & 0 & 1 \\ 0 & 0 & 0 \\ 1 & 0 & 0 \end{array} \right]  V'_3 =
\frac{1}{3(1 + \sqrt{3})} \left[ \begin{array}{cc} 2 & -(1+\sqrt{3}) \\ -(1+\sqrt{3}) & -4-2\sqrt{3}  \end{array} \right],
\]
which has eigenvalues $-1$ and $1/3$.
Therefore, it follows from  the interlacing theorem that $\mmax \geq 1/3$ and thus ($-V^TAV$) is not PSD.
Consequently, $A$ is not an EDM.
\epr

A remark is in order here. Let $G$ be a complete multipartite graph and assume that its nodes are
partitioned into $s$ independent sets. Then
the adjacency matrix of $G$ is an EDM embedding dimension $s-1$. For example, the nodes of $K_n$ and $\overline{K_n}$
are, obviously, partitioned into $n$ and $1$ independent sets. Consequently, the embedding dimensions of the corresponding adjacency
matrices, i.e., $E-I$ and $\bz$, are respectively $n-1$ and $0$ as expected.

\begin{thm} \label{thmmin-1Kn}
Let $G$ be a graph on $n$ nodes which is not null and let  $A$ denote its adjacency matrix.
Then  $\mmin \leq -1$. Moreover, $\mmin = -1$ if and only if
$G$ is a cluster graph.
\end{thm}

\bpr Assume that $G$ is a cluster graph.
If $G$ is the disjoint union of $K_2$ and ($n-2$) isolated nodes, then by the proof of Theorem \ref{thmcmpar}, we have  $\mmin = -1$.
Otherwise, it follows from Proposition \ref{propinterl} that $\mmin = -1$ since $(-1)$ is an eigenvalue of $A$ of
multiplicity at least 2.

To prove the other direction, assume that $G$ is not a cluster graph. Then
$G$ has $P_3$ as an induced subgraph.
Therefore, wlog, assume that the nodes of $P_3$ are 1, 2, and 3. Therefore, the  third leading principal submatrix of $A$ is
\[
 \left[ \begin{array}{ccc} 0 & 1 & 1 \\ 1 & 0 & 0 \\ 1 & 0 & 0 \end{array} \right].
\]
Let  $V$ be as defined in (\ref{defVs}). Then
the second leading  principal submatrix of $V^TAV$ is
\[
{V'}_3^T  \left[ \begin{array}{ccc} 0 & 1 & 1 \\ 1 & 0 & 0 \\ 1 & 0 & 0 \end{array} \right]  V'_3 =
\frac{-2}{3} \left[ \begin{array}{cc} 1 & 1 \\ 1 & 1  \end{array} \right].
\]
Therefore, it follows from the interlacing theorem  that $\mmin \leq - 4/3$.
\epr

Let $G$ be a cluster graph and assume that $G \neq K_n$ and $G \neq  \overline{K_n}$. Then
$G$ is the disjoint union of at least two complete graphs, say $K_{n_1}$ and $K_{n_2}$, where
$n_1 \geq 2$. Thus $G$ has an induced $\overline{P_3}$ whose nodes are  two  from $K_{n_1}$ and one
 from $K_{n_2}$. Consequently, $G$ is not a complete multipartite graph. On the other hand,
if $G = K_n$, then obviously $\mmax=\mmin=-1$ and if
$G = \overline{K_n}$, then trivially $\mmax=\mmin=0$.
Hence, we have proven the following corollary.

\begin{cor} \label{cornoopenends}
There exists no graph $G$ such that $\mmax=0$ and $\mmin = -1$.
\end{cor}

\section{Euclidean Representations}

let $G$ be a simple graph on $n$ nodes which is neither complete nor null.
Then $G$ admits a Euclidean representation in $\Rs^r$ iff there exist two distinct positive scalars
$\alpha$ and $\beta$ such that $D = \alpha A + \beta \bar{A}$
is an EDM of embedding dimension $r$.
Wlog assume that $\alpha = 1$ and thus $0 < \beta \neq 1$.
Hence, $D = \beta(E-I) + (1 - \beta) A$.
Next, we derive upper and lower bounds on $\beta$ such that $D$ is an EDM. To this end,
 $X$, the projected Gram matrix of $D$, is given by
\beq \label{XbetaA}
2 X = \beta I_{n-1} + (\beta - 1) V^T A V.
\eeq
Hence, by Theorem \ref{thmbasic}, $D$ is an EDM of dimension $r$ iff $X$ is PSD of rank $r$.
Recall that $\mmin$ and $\mmax$ denote, respectively, the minimum and the maximum eigenvalues of $V^TAV$.
Assume that $\beta > 1$. Then, in light of Theorem \ref{thmmin-1Kn}, $X$ is PSD iff
\begin{eqnarray*}
1 < \beta <   + \infty && \mbox{if  $\mmin = -1$}, \\
    1 < \beta \leq  \frac{|\mmin |}{|\mmin | - 1}   && \mbox{if  $\mmin < -1$}.
\end{eqnarray*}
On the other hand, assume that $0 < \beta < 1$. Then, in light of Lemma \ref{lemmaxE-I} and since $G \neq K_n$, $X$ is PSD iff
\begin{eqnarray*}
 0 < \beta < 1  && \mbox{if  $\mmax = 0$}, \\
\frac{\mmax }{\mmax  + 1} \leq \beta < 1  && \mbox{if  $\mmax > 0$}.
\end{eqnarray*}

Let us define
\beq \label{defbetaul}
\beta_l = \frac{\mmax }{\mmax  + 1}  \mbox{ and } \beta_u = \frac{|\mmin |}{|\mmin | - 1}.
\eeq
Therefore, $X$ is PSD iff
\beq
\begin{array}{ll}
  \beta \in [\beta_l, 1 ) \cup (1, + \infty) & \mbox{ if } \mmin = -1 \mbox{ and } \mmax > 0  , \\
  \beta \in (0, 1) \cup (1, \beta_u] & \mbox{ if } \mmin < -1 \mbox{ and } \mmax = 0  , \\
  \beta \in [\beta_l, 1) \cup (1, \beta_u] & \mbox{ if } \mmin < -1 \mbox{ and } \mmax > 0.
\end{array}
\eeq
Note that Corollary \ref{cornoopenends} rules out the case in which $\mmin = -1$ and $\mmax=0$.
Therefore, for $G \neq K_n$ and $G \neq \overline{K_n}$, we have
\beq
\begin{array}{ll}
  \beta \in [\beta_l, 1 ) \cup (1, + \infty) & \mbox{ if $G$ is a cluster graph}, \\
  \beta \in (0, 1) \cup (1, \beta_u] & \mbox{ if $G$ is a complete multipartite graph}, \\
  \beta \in [\beta_l, 1) \cup (1, \beta_u] & \mbox{ otherwise}.
\end{array}
\eeq

Let $m(\mmin)$ and $m(\mmax)$ denote the multiplicities of $\mmin$ and $\mmax$. Then
\beq \label{eqrankX}
\rank(X) = \left\{ \begin{array}{ll} n-1 & \mbox{ if } \beta \neq \beta_l \mbox{ and } \beta \neq \beta_u, \\
                                     n-1-m(\mmax)  & \mbox{ if } \beta = \beta_l, \\
                                      n-1-m(\mmin)  & \mbox{ if } \beta = \beta_u.
                                     \end{array} \right.
\eeq
Therefore, if  $G$ is a cluster graph, then $\dE = n-1-m(\mmax)$; and if $G$ is a complete multipartite graph,
then $\dE = n-1-m(\mmin)$. Otherwise,

\beq \label{dimE1st}
\dE =   \min \{n-1- m(\mmax) ,n-1- m(\mmin) \}.
\eeq
As a result, $n-2$ is an upper bound on $\dE$ as proved in \cite{ngu18}.
Using a different approach, Einhorn and Schoenberg \cite{es66a, es66b}, obtained an equivalent equation for $\dE$.
Next, we derive their equation and we show the equivalence between the two equations. To this end,
let $D(t)= A + t \bar{A}$ and let
\[
\tilde{X}(t) = - [ -e_{n-1} \;\; I_{n-1}] D(t)  \left[ \begin{array}{c} -e_{n-1}^T \\ I_{n-1} \end{array} \right].
\]
Now $\left[ \begin{array}{c} -e_{n-1}^T \\ I_{n-1} \end{array} \right] = V \Phi$ for some nonsingular $\Phi$ since
the columns of $\left[ \begin{array}{c} -e_{n-1}^T \\ I_{n-1} \end{array} \right]$ form a basis of $e^{\perp}$.
Therefore,
$\tilde{X} (t) =2 \Phi^T X(t) \Phi$, where $X(t)$ is the projected Gram matrix of $D(t)$.
Thus, $\tilde{X}(t)$ is PSD and of rank $r$ iff $X(t)$ is PSD and of rank $r$.
The {\em discriminating polynomial} \cite{es66a} of $D(t)$ is defined as
\beq  \label{dispoly}
p(t) = \det(\tilde{X}(t)).
\eeq
Note that $\tilde{X}(1) = I_{n-1} + E_{n-1}$ is PD and thus $p(1) > 0$.
Also, note that $D(t)$ is an EDM iff $\tilde{X}(t)$ is PSD.
Now if all roots of $p(t)$ are $< 1$, then $\tilde{X}(t)$ is PSD for all $t \geq 1$.
Otherwise, $\tilde{X}(t)$ is PSD for all $t : 1 \leq t \leq t_2$, where
$t_2$ is the smallest root of $p(t)$ such that $t_2 > 1$.
On the other hand, let $t_1$ be the largest root of $p(t)$ such that $0< t_1 <1$ if such root exists.
Thus, $\tilde{X}(t)$ is PSD for all $t_1 \leq t \leq 1$ if $t_1$ exists, and
for all $0 < t \leq 1$  otherwise.
As a result,
Einhorn and Schoenberg obtained that, if both $t_1$ and $t_2$ exist, then
\beq \label{dimEpoly}
\dE =  \min \{n-1- m(t_1) ,n-1- m(t_2) \},
\eeq
where $m(t_1)$ and $m(t_2)$ are the multiplicities of $t_1$ and $t_2$.
To establish the equivalence between Equations (\ref{dimE1st}) and (\ref{dimEpoly}),
note that
\[
p(t) = c \det (t I_{n-1} + (t - 1) V^T A V) = c (t-1)^{n-1} \chi( \frac{t}{1-t}),
\]
where $c$ is a constant and  $\chi(\mu)$ is the characteristic polynomial of $V^TAV$.
Hence,
\[
t_1  = \frac{\mmax}{\mmax + 1} \mbox{ and } t_2 = \frac{|\mmin|}{|\mmin| - 1}.
\]

Let $G$ be a $k$-regular graph, which is neither complete nor null, and
let $\lambda_1(A)=k \geq \lambda_2(A) \geq \cdots$ $ \geq \lambda_n(A)$
be the eigenvalues of its adjacency matrix. Then Proposition \ref{propreg} implies that $\mmin=\lambda_n(A)$ and $\mmax= \lambda_2(A)$.
Hence, as was proven in \cite{roy10}, if $G$  is neither a cluster graph nor a complete multipartite graph, then
\beq
\dE = \min \{ n-1- m(\lambda_2(A)) ,n-1- m(\lambda_n(A)) \}.
\eeq
\begin{exa} \label{exC5}
Let $G=C_5$, the cycle on 5 nodes. Then $\mmax = (\sqrt{5} - 1)/2$ with multiplicity 2 and
$\mmin = -(\sqrt{5} + 1)/2$ with multiplicity 2. Hence, $\dE = 2$,  $\beta_u = (\sqrt{5}+3)/2$ and
$\beta_l = (-\sqrt{5}+3)/2$. Observe that $\beta_u \beta_l = 1$ as expected since the regular pentagon is
the unique two-distance representation of $C_5$ in $\Rs^2$.
\end{exa}

Let $D= A + \beta \bar{A}$ be the EDM of a Euclidean representation of $G$.
Obviously, every Euclidean representation of $G$ is at the same time a Euclidean representation of the complement graph $\bar{G}$.
More precisely,
\beq
   \bar{D} =  \frac{D}{\beta} = \bar{A} +  \frac{1}{\beta} A
\eeq
is the EDM of a Euclidean representation of $\bar{G}$. Moreover, if we let $\bar{\mu}_{\min}$ and $\bar{\mu}_{\max}$
denote, respectively, the minimum  and the maximum eigenvalues of $V^T \bar{A}V$.
Then it is immediate that
\beq \label{minmaxAAb}
 \bar{\mu}_{\min} = -1 - \mmax \;\; \mbox{ and } \;\; \bar{\mu}_{\max} = -1 - \mmin,
\eeq
and
\beq  \label{minmaxAAbm}
m (\bar{\mu}_{\min})= m(\mmax) \;\; \mbox{ and } \;\; m(\bar{\mu}_{\max}) = m(\mmin).
\eeq
Consequently, $\dE$ = dim$_E(\bar{G})$ as expected.
We end this section by noting the following well-known lower bound on $\dE$.
It is well known \cite{blo84,bbs83} that any two-distance $n$-point configuration in $\Rs^r$ satisfies
\[
 n \leq \frac{(r+1)(r+2)}{2}.
\]
Hence, for any graph $G$, which is neither complete nor null, we have
\[
\dE \geq \frac{1}{2} (\sqrt{8n+1} - 3).
\]

\section{Spherical Representations}

Graph $G$ admits a spherical representation in $\Rs^r$ iff there exist two distinct positive scalars
$\alpha$ and $\beta$ such that $D = \alpha A + \beta \bar{A}$ is a spherical EDM of embedding dimension $r$.
Wlog, assume that $\alpha=1$ and hence $D= (1-\beta)A + \beta(E-I)$.
Assume that $D$ is an EDM and let  $G$ be a $k$-regular graph, i.e., $Ae=ke$.
Then, $De= ((1-\beta)k + \beta(n-1)) e$. Consequently, $D$ is a regular EDM with radius
\[
\rho = \left(\frac{(1-\beta)k+\beta(n-1) }{2n} \right)^{1/2}.
\]

Now for a general graph $G$, Theorem \ref{thmspher} and equation (\ref{eqrankX}) imply that
if $\beta \neq \beta_l$ and $\beta \neq \beta_u$, where $\beta_l$ and $\beta_u$ are as defined in (\ref{defbetaul}),
then $D$ is spherical of radius
\[
\rho = \left(\frac{1} {2 e^T((1-\beta)A +\beta(E-I))^{-1}e } \right)^{1/2},
\]
since, in this case,  $w=D^{-1}e$.
Otherwise, if
$\beta = \beta_u$ or $\beta = \beta_l$, then $D$ may or may not be spherical.
To characterize the sphericity of $D$ in this case, we need the following two definitions.

\begin{defn} \label{defnUu}
 For adjacency matrix $A$, let $U_u$ be the $(n-1) \times m(\mmin))$ matrix whose columns form an orthonormal basis
for the eigenspace of $V^TAV$ associated with $\mmin$. That is, the columns of $U_u$ are orthonormal eigenvectors
of $V^TAV$ corresponding to $\mmin$.
\end{defn}

\begin{defn} \label{defnUl}
 For adjacency matrix $A$, let $U_l$ be the $(n-1) \times m(\mmax))$ matrix whose columns form an orthonormal basis
for the eigenspace of $V^TAV$ associated with $\mmax$.
\end{defn}

The following theorem establishes a necessary and sufficient condition for EDMs $D_u = A + \beta_u \bar{A}$ and
 $D_l = A + \beta_l \bar{A}$ to be spherical.

\begin{thm}
Let $D_u = A + \beta_u \bar{A}$, where $\beta_u$ is as given in (\ref{defbetaul}).
Then the EDM $D_u$ is spherical if and only if
\[
  A VU_u = \mmin VU_u.
\]
Similarly, let $D_l = A + \beta_l \bar{A}$.
Then the EDM $D_l$ is spherical if and only if
\[
  A VU_l = \mmax VU_l.
\]
\end{thm}

\bpr
We present the proof for $D_u$. The proof for $D_l$ is similar.
Now it follows from (\ref{XbetaA}) that
the null space of $X_u$, the projected Gram matrix of $D_u$, is given by
\[
\mbox{null} (X_u)= \left\{ \xi \in \Rs^{n-1}: V^T A V \xi = \frac{\beta_u}{1-\beta_u} \xi = \mmin \xi \right\}.
\]
Moreover, it also follows from (\ref{XbetaA}) that
\[
 2 X_u U_u = ( (1+ \mmin) \beta_u  - \mmin) U_u = \bz.
\]
Hence, the columns of $U_u$ form an orthonormal basis of the null space of $X_u$ and thus, by Lemma \ref{lemZVU},
$VU_u$ is a Gale matrix of $D_u$. Therefore, by Theorem \ref{thmspher}, $D_u$ is spherical iff $D_u V U_u = \bz$; i.e.,
iff
\[
 A VU_u = \frac{\beta_u}{1 - \beta_u} VU_u = \mmin VU_u.
\]
\epr

The following corollaries are immediate.

\begin{cor} \label{corclusterS}
Let $G$ be a graph on $n$ nodes, $G \neq K_n$ and $G \neq \overline{K_n}$. Let $Z_l = VU_l$.
If $G$ is a cluster graph, i.e., if $\mmin = -1$, then
\[
\dS = \left\{ \begin{array}{ll} n-1-m(\mmax) & \mbox{ if } A Z_l = \mmax Z_l, \\
                               n-1 & \mbox{ otherwise}. \end{array} \right.
\]
\end{cor}

\begin{cor}  \label{corpartiteS}
Let $G$ be a graph on $n$ nodes, $G \neq K_n$ and $G \neq \overline{K_n}$. Let $Z_u = VU_u$.
If $G$ is a complete multipartite graph, i.e., if $\mmax = 0$, then
\[
\dS = \left\{ \begin{array}{ll} n-1-m(\mmin) & \mbox{ if } A Z_u = \mmin Z_u, \\
                               n-1 & \mbox{ otherwise}. \end{array} \right.
\]
\end{cor}

We should point out that, in light of (\ref{minmaxAAbm}),
 $m(\mmin) $  in Corollary \ref{corclusterS} is equal to $m(\mmax)$ in Corollary \ref{corpartiteS}.
Furthermore, it is easy to see that
$Z_l$  in Corollary \ref{corclusterS} is equal to $Z_u$ in Corollary \ref{corpartiteS},
and that
$A Z_l = \mmax Z_l $  in Corollary \ref{corclusterS} iff  $ A Z_u = \mmin Z_u$ in Corollary \ref{corpartiteS}.

\begin{cor}
Let $G$ be a graph on $n$ nodes, $G \neq K_n$ and $G \neq \bar{K_n}$.
Assume that $G$ is neither a cluster graph nor a complete multipartite graph, i.e.,  $\mmin \neq -1$ and  $\mmax \neq 0$.
Let $Z_l = VU_l$ and $Z_u = VU_u$. Then
\[
\dS = \left\{ \begin{array}{ll} n-1 & \mbox{ if } A Z_u \neq \mmin Z_u \mbox{ and } A Z_l \neq \mmax Z_l, \\
                               n-1-m(\mmin) & \mbox{ if } A Z_u = \mmin Z_u \mbox{ and } A Z_l \neq \mmax Z_l, \\
                               n-1-m(\mmax) & \mbox{ if } A Z_u \neq \mmin Z_u \mbox{ and } A Z_l = \mmax Z_l.
                               \end{array} \right.
\]
Otherwise, i.e.,  if $A Z_u = \mmin Z_u \mbox{ and } A Z_l = \mmax Z_l$. Then
\[
\dS = \min \{ n-1-m(\mmin), n-1-m(\mmax)\}.
\]
\end{cor}

\begin{exa}
Let $G=C_5$ which was considered in Example \ref{exC5}, where  $\mmax = (\sqrt{5} - 1)/2$ with multiplicity 2 and
$\mmin = -(\sqrt{5} + 1)/2$ with multiplicity 2. Now since $G$ is $2$-regular, $D = A + \beta \bar{A}$ is a spherical
EDM for all $\beta_l \leq \beta \leq \beta_u$. On the other hand, two Gale matrices of  
$D_l = A + \beta_l \bar{A}$ and $D_u=A + \beta_u \bar{A}$ are
\[
 Z_l  = \left[ \begin{array}{cc} 1 & 0 \\ 0 &  1 \\ -1 & \mmax \\ -\mmax & -\mmax \\ \mmax & -1 \end{array} \right]  \;\; \mbox{ and } \;\;
 Z_u  = \left[ \begin{array}{cc} 1 & 0 \\ 0 & 1 \\ -1 & \mmin \\ -\mmin & -\mmin \\ \mmin & -1 \end{array} \right].
\] 
It is easy to verify that $A Z_l = \mmax Z_l$ and $A Z_u = \mmin Z_u$.
Therefore, as expected, both EDMs $D_l$ and $D_u$ are spherical with radii (squared)
$\rho_l^2 = 2/(5+\sqrt{5})$, and $\rho_u^2 = 2/(5-\sqrt{5})$.
As a result, $\dS = 2$.
\end{exa}

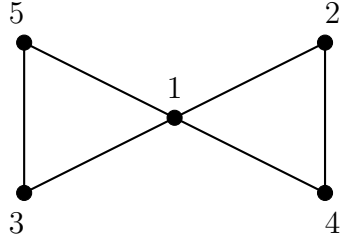
\begin{figure}[t!]
\begin{center}
\tikzstyle{place}=[circle,draw=black!100,fill=black!100,thick]
\tikzstyle{transition}=[rectangle,draw=black!50,fill=black!20,thick]
\tikzstyle{blank}=[circle,draw=black!100,fill=black!100,thick]
\begin{tikzpicture}[auto, inner sep=.5mm]
\draw[fill=black!100] (0,0) circle(0.1 cm);
\draw[fill=black!100] (2,1) circle(0.1 cm);
\draw[fill=black!100] (-2,-1) circle(0.1 cm);
\draw[fill=black!100] (2,-1) circle(0.1 cm);
\draw[fill=black!100] (-2,1) circle(0.1 cm);

\draw[thick] (2,1)--(2,-1);
\draw[thick] (-2,1)--(-2,-1);
\draw[thick] (0,0)--(2,-1);
\draw[thick] (0,0)--(2,1);
\draw[thick] (0,0)--(-2,-1);
\draw[thick] (0,0)--(-2,1);

\draw (0,0.4) node {$1$};
\draw (2.1,1.4) node {$2$};
\draw (-2.1,-1.4) node {$3$};
\draw (2.1,-1.4) node {$4$};
\draw (-2.1,1.4) node {$5$};

\end{tikzpicture}
\caption{The bow tie graph of Example \ref{exsphernonspher}.}
\label{figsphernonspher}
\end{center}
\end{figure}

\begin{exa} \label{exsphernonspher}
Consider the ``bow tie" graph depicted in Figure \ref{figsphernonspher}. Then, in this case, $\mmin=-1.4$, $\mmax=1$ and
$m(\mmin) = m(\mmax) = 1$.
Thus, $\beta_l=1/2$ and $\beta_u = 7/2$. Moreover,
\[
 Z_l = V U_l = \frac{1}{2} \left[ \begin{array}{r} 0 \\ 1 \\ -1 \\ 1 \\ -1 \end{array} \right] \;\; \mbox{ and } \;\;
 Z_u = V U_u = \frac{1}{\sqrt{20}} \left[ \begin{array}{r} -4 \\ 1 \\ 1 \\ 1 \\ 1 \end{array} \right].
\]
It is easy to verify that $A Z_l = \mmax Z_l$ and $A Z_u \neq \mmin Z_u$.
Therefore, the EDM $D_l = A + \beta_l \bar{A}$ is spherical of radius $\rho_l = 1/\sqrt{3}$, while
the EDM $D_u=A + \beta_u \bar{A}$ is not spherical.
As a result, $\dS = 3$.
\end{exa}

Obviously, every spherical representation of $G$ is at the same time a spherical representation of the complement graph $\bar{G}$.
More precisely,
\[
   \bar{D} =\frac{D}{\beta} = \bar{A} +  \frac{1}{\beta} A
\]
is the EDM of a spherical representation of $\bar{G}$. Therefore, it follows from (\ref{minmaxAAb}) and (\ref{minmaxAAbm}) that
$\dS$ = dim$_S(\bar{G})$ as expected. Moreover, if we let $\bar{\rho}$ denote the radius of $\bar{D}$, then clearly
$\bar{\rho}^2 = \rho^2 / \beta$.

Now suppose that the EDM $D_u = A + \beta_u \bar{A}$ is spherical. Then $\rho$, the radius of the sphere containing its generating points,
can be given explicitly in terms of $A$. To this end, let $P$ be a configuration matrix of $D_u$ such that $P^T e = \bz$ and let
$B_u=PP^T$.
By Theorem \ref{thmspher}, $\rho^2 = a^Ta + e^T D_u e / (2 n^2)$, where $2Pa = (I - E/n)  \diag(B_u)$.
Now it follows from (\ref{defs}), since $s=e/n$, that
\[
\diag(B_u) = \frac{D_ue}{n} - \frac{e^TD_ue}{2 n^2} e.
\]
Therefore,
\[
4 a^T a = (\diag(B_u))^T P (P^TP)^{-2} P^T \diag(B_u) = \frac{e^T D_u B_u^{\dag} D_u e}{n^2},
\]
  where $B_u^{\dag} =  P (P^TP)^{-2} P^T$ is the Moore-Penrose inverse of $B_u$.
Let $V^T A V = \mmin U_u U_u^T + W_u \Lambda_u W_u^T$ be the spectral decomposition of $V^TAV$.
Then the projected Gram matrix of $D_u$ is given by
\[
 X_u = \frac{1}{2 (\mmin + 1) } W_u ( \mmin I - \Lambda_u) W_u^T,
\]
and thus
\[
X_u^{\dag} = 2 (\mmin + 1) W_u (\mmin I - \Lambda_u)^{-1} W_u^T.
\]
Hence,
\[
 a^T a = \frac{1}{2 n^2 (\mmin + 1) } e^T A V W_u (\mmin I -\Lambda_u)^{-1} W_u^T V^T A e.
\]
As a result,
\[
\rho^2 = \frac{1}{2 n^2 (\mmin + 1) } ( e^T A V W_u (\mmin I -\Lambda_u)^{-1} W_u^T V^T A e +  \mmin(n^2-n) + e^T Ae).
\]

We end this section by noting the following well-known lower bound on $\dS$.
It is well known \cite{dgs77} that any two-distance $n$ point spherical configuration in $\Rs^r$ satisfies
\[
 n \leq \frac{r(r+3)}{2}.
\]
Hence, for any graph $G$ we have
\[
\dS \geq \frac{1}{2} (\sqrt{8n+9} - 3).
\]

\section{J-Spherical Representations}

Musin \cite{mus18} proved that every graph $G$, which is not complete or null, admits a unique, up to an isometry,
$J$-spherical representation. Unfortunately, his proof is not constructive.
In this section, we give explicit simple formulas for the $J$-spherical representation of $G$
and for  $\dJ$ in terms of the largest eigenvalue of $\bar{A}$, the adjacency matrix of the complement graph $\bar{G}$,
and its multiplicity. We also answer a question raised in \cite{mus18}.

Evidently, $G$ admits a $J$-spherical representation in $\Rs^r$ iff there exists
a scalar $\beta > 2$ such that $D = 2 A + \beta \bar{A}$ is a spherical EDM of unit radius and of embedding dimension $r$.
Let $D$ be a spherical EDM of unit radius. Then Theorem \ref{thmspher} implies that
$2 e^Tw = 1$ where $Dw=e$. Consequently, we will find it convenient, in this section,  to set $s=2 w$ in Theorem \ref{thmbasic};
i.e., we fix the origin such that $B$, the Gram matrix of $D$, satisfies $Bw=\bz$.
Therefore, $B$, in this case, is given by
\beq
 B= E - \frac{1}{2}D \;\; \mbox{ and satisfies } Bw = \bz.
\eeq
Let $\beta = 2 + 2\delta$, where $\delta > 0$. Then $D= 2 A + \beta \bar{A} = 2 (E-I) + 2 \delta \bar{A}$ and thus
\beq
 B = I  -  \delta \bar{A} \mbox{ and satisfies } Bw = \bz.
\eeq
As a result, $G$ admits a $J$-spherical representation in $\Rs^r$ iff there exists
$\delta > 0 $ such that
\beq
B= I  -  \delta \bar{A} \mbox{ is PSD }, \;\; Bw = \bz \mbox{ and } \rank(B)=r.
\eeq

Now let $\lambda_1(\bar{A}) \geq \cdots \geq \lambda_n(\bar{A})$ be the eigenvalues of $\bar{A}$.
Then $B$ is PSD iff
\[
\delta \leq \frac{1}{\lambda_1(\bar{A})}.
\]
On the other hand, $Bw= \bz$ is equivalent to
\[
  \bar{A} w = \frac{1}{\delta} w.
\]
Hence, $1/ \delta$ is an eigenvalue of $\bar{A}$ and thus $1/\delta \leq \lambda_1(\bar{A})$.
Consequently,
\beq
\delta = \lambda_1(\bar{A}).
\eeq
As a result, we have proven the following theorem.
\begin{thm} \label{thmmainJ}
Let $G$ be a graph on $n$ nodes, which is neither complete nor null, and let $\delta = 1/\lambda_1(\bar{A})$, where
$\lambda_1(\bar{A})$ is the largest eigenvalue of $\bar{A}$, the adjacency matrix of the complement graph $\bar{G}$.
Then $G$ admits a unique, up to an isometry, $J$-spherical representation whose EDM is given by
\[
 D = 2(E - I) + 2 \delta \bar{A}.
\]
Moreover, $\dJ= n - m(\lambda_1(\bar{A}))$.
\end{thm}

Following \cite{mus18}, let us refer to $\alpha = 2$ as the first  distance (squared) and to
$\beta = 2 + 2 \delta$ as the second distance (squared).
The following observation is worth pointing out.
It follows from the interlacing theorem and (\ref{minmaxAAb})
that $\lambda_1(\bar{A}) \geq \mmaxb =  |\mmin| - 1$. Hence,
the second distance (squared) satisfies
\[
\beta \leq 2 \frac{|\mmin|}{|\mmin| - 1} = 2 \beta_u
\]
as expected. Note that the factor of $2$ results from the fact that $\alpha$,
the first distance (squared), is 2 instead of 1 as was the case in previous sections.

\begin{exa}
Consider the graph $G= C_5$. Then $\lambda_1(\bar{A}) = 2 $ with multiplicity $1$.
Thus, as was observed in \cite{mus18}, $\dJ = 4$.
\end{exa}

\begin{exa}
Let $G$ be the ``bow tie" graph depicted in Figure \ref{figsphernonspher} and considered in Example \ref{exsphernonspher}.
Then $\lambda_1(\bar{A}) = 2 $ with multiplicity $1$. Thus, $\delta = 1/2$ and $\dJ = 4$.
\end{exa}

\begin{exa}
Let $G = K_{n_1,\ldots,n_s}$ be a complete multipartite graph, where $n=n_1+\cdots+n_s$ and
\[
n_1 = \cdots = n_k > n_{k+1} \geq \cdots \geq n_s.
\]
Then Musin \cite{mus18} proved that $\dJ = n - k$. Clearly, $\bar{G}$ in this case is a cluster graph and
$\lambda_1(\bar{A})= n_1 - 1$ with multiplicity $k$.
\end{exa}

We conclude this paper by presenting a characterization of graphs whose $J$-spherical representations
 have the same second distance (squared) $\beta =2+2 \delta$.
This characterization follows as an immediate corollary of Theorem \ref{thmmainJ} and answers a question raised by Musin
\cite{mus18}.

\begin{thm} \label{thmsamebeta}
Let $G_1$ and $G_2$ be two distinct graphs, which are neither complete nor null, and let $\bar{A_1}$ and
$\bar{A_2}$ be, respectively, the adjacency matrices of the complement graphs $\bar{G_1}$ and $\bar{G_2}$.
Then the two $J$-spherical representations of $G_1$ and $G_2$ have the same second distance (squared) if and only if
 $\lambda_1(\bar{A_1}) = \lambda_1(\bar{A_2})$.
\end{thm}

We conclude this paper with the following two examples as an illustration of Theorem \ref{thmsamebeta}.

\begin{exa}
Musin \cite{mus18} gave the following cluster graphs
$G_1= 3K_2$, $G_2 = 2K_4$, $G_3 = K_2 \cup K_8$ and $G_4= K_1 \cup K_{16}$ as an example of graphs whose
$J$-spherical representations have the same second (squared) distance of $\beta = 5/2$. It is easy to verify
that for all these graphs $\lambda_1(\bar{A}) = 4$ and thus $\beta = 5/2$.
\end{exa}

\begin{exa}
Consider the graphs $G_n = \overline{C_n}$ for $n \geq 4$. It is immediate that for all these graphs $\lambda_1(\bar{A}) = 2$.
Hence, the $J$-spherical representations of all these graphs have the same second distance (squared) of $\beta = 3$.
\end{exa}



\end{document}